\documentclass{svjour2}                    % onecolumn
\smartqed 

\usepackage{hanswill}

\begin{document}

\title{On the Hopf Algebraic Structure of\\ Lie Group
  Integrators}

\journalname{Foundations of Computational Mathematics}

\titlerunning{Hopf algebra of ordered trees}      

\author{H.~Z.~Munthe-Kaas \and W.~M.~Wright}

\authorrunning{Munthe-Kaas \and Wright} 

\institute{H.~Z.~Munthe-Kaas \at
              Department of Mathematics, University of Bergen,
              N-5020 Norway\\
              \email{H.Munthe-Kaas@mi.uib.no}
           \and
           W.~M.~Wright \at
              Department of Mathematics, La Trobe University,
              Victoria 3086, Australia\\
              \email{W.Wright@latrobe.edu.au}
}

\date{Received: date / Accepted: date}

\maketitle

\begin{abstract}
  A commutative but not cocommutative graded Hopf algebra $\Hn$, based
  on ordered rooted trees, is studied. This Hopf algebra generalizes
  the Hopf algebraic structure of unordered rooted trees $\Hc$,
  developed by Butcher in his study of Runge--Kutta methods and later
  rediscovered by Connes and Moscovici in the context of
  non-commutative geometry and by Kreimer where it is used to describe
  renormalization in quantum field theory. It is shown that $\Hn$ is
  naturally obtained from a universal object in a category of
  non-commutative derivations, and in particular, it forms a
  foundation for the study of numerical integrators based on
  non-commutative Lie group actions on a manifold.  Recursive and
  non-recursive definitions of the coproduct and the antipode are
  derived. It is also shown that the dual of $\Hn$ is a Hopf algebra
  of Grossman and Larson.  $\Hn$ contains two well-known Hopf algebras
  as special cases: The Hopf algebra $\Hc$ of Butcher--Connes--Kreimer
  is identified as a proper subalgebra of $\Hn$ using the image of a
  tree symmetrization operator.  The Hopf algebra $\Hf$ of the Free
  Associative Algebra is obtained from $\Hn$ by a quotient
  construction.

\end{abstract}

% ===========================================================================
\section{Introduction} 
\label{sect:intro}
% ===========================================================================

Since Cayley~\cite{cayley57ott} in 1857 and 100 years later Merson
\cite{merson57aom} and followed shortly after by Butcher
\cite{butcher63cft}, it has been known that rooted trees are extremely
useful for structuring algebras of differential operators and
elementary differentials. In 1972 Butcher~\cite{butcher72aat} produced
the far-reaching result that Runge--Kutta methods form a group. This
was later named the Butcher group in the paper by Hairer and Wanner
\cite{hairer74otb}, who made significant contributions to this
theory. The Butcher group is defined on the dual of the tree space and
it was pointed out by D\"ur~\cite{dur86mfi} that there exists a
one-to-one correspondence between the Butcher group and the
commutative graded Hopf algebra of unordered rooted trees. The Hopf
algebra of unordered rooted trees has had far-reaching applications in
various areas of mathematics and physics. In 1998 a Hopf subalgebra
was discovered by Connes and Moscovici~\cite{connes98hac} during work
on an index theorem in non-commutative geometry and by
Kreimer~\cite{kreimer98oth} in the renormalization method of quantum
field theory. Further collaborations between Connes and
Kreimer~\cite{connes98har,connes00riq} have lead to other surprising
results; notably a connection with the Riemann--Hilbert problem was
established. Brouder~\cite{brouder00rkm,brouder04tra} realized that
the mathematical structure of Connes--Kreimer was the same as that of
Butcher. Grossman and Larson~\cite{grossman89has} also developed a
cocommutative graded Hopf algebra on a general class of rooted
trees. It was shown by Foissy~\cite{foissy02ladII} and Hoffman
\cite{hoffman03cor} that the commutative Hopf algebra of Butcher and
Connes--Kreimer was the dual of the cocommutative Hopf algebra of
Grossman--Larson, which corrected the original result of
Panaite~\cite{panaite00rtc}. Murua \cite{murua99fsaI,murua05tha} has
developed series expansions of elementary differential operators and
shown among other results that the logarithm of such a series is
equivalent to the series expansions obtained from backward error
analysis.

Recently a great deal of interest has been focused on developing
numerical methods which preserve geometric properties of the exact
flow. In particular, Lie group integrators, which describe integrators
that use Lie group actions on manifolds, were originally proposed by
Crouch and Grossman \cite{crouch93nio} followed shortly after by Lewis
and Simo \cite{lewis94caf,lewis95caf}. Integrators of this type are
now known as Lie group integrators; a survey of these methods is given
in \cite{iserles00lgm}. Series expansions for various classes of Lie
group methods have been developed; these expansions are generally used
to analyze order. Munthe-Kaas \cite{munthe-kaas95lbt,munthe-kaas98rkm}
constructed the order conditions for a special subclass of Lie group
methods, where the computations are performed in a Lie algebra, which
is a linear space. Later, \cite{munthe-kaas99hor} it was shown that
the classical order conditions could be used along with a certain
transformation. Owren and Marthinsen \cite{owren99rkm} developed the
general order conditions for the Crouch--Grossman methods with their
analysis being based on ordered rooted trees. Recently, Owren
\cite{owren05ocf} derived the order conditions for the commutator free
Lie group methods \cite{celledoni03cfl}, which were derived to
overcome some of the problems associated with computing commutators.

In this paper we aim to construct a commutative graded Hopf algebraic
structure, which can be used to analyze the order of all Lie group
methods. The outline of this paper is as follows: In
Section~\ref{sect:derivation} we will introduce ordered trees and
forests and describe some useful operations on them. We will motivate
the present Hopf algebra as a universal object in a category of
non-commutative derivations, and also briefly discuss Lie--Butcher
theory which will be treated in more detail in
\cite{munthe-kaas05lbt}. In Section~\ref{sect:hopf} we develop the
Hopf algebra of ordered trees, giving both recursive and non-recursive
definitions of the coproduct and antipode, using certain cutting
operations on ordered rooted trees.  We show that the Hopf algebra
described in this paper is the dual of the Grossman--Larson Hopf
algebra, thus generalizing the result of Hoffman~\cite{hoffman03cor}.
Finally, in Section~\ref{sect:homo} we use a symmeterization operator
to provide an injective Hopf algebra homomorphism from the unordered
trees into the ordered trees, establishing the former as a sub-algebra
of the latter.

% ===========================================================================
\section{Algebras of non-commutative derivations}
\label{sect:derivation}
% ===========================================================================

% ===========================================================================
\subsection{An algebra of trees}
\label{sect:algebra}
% ===========================================================================

In this section we will define an algebra $N$ spanned by forests of
ordered (and possibly colored) rooted trees. This algebra is a
universal (`free') object in a general category of non-commutative
derivation algebras, and plays a role in symbolic computing with
Lie--Butcher series similar to the role of \emph{free Lie
  algebras}~\cite{munthe-kaas99cia,reutenauer93fla} in symbolic
computing with Lie algebras.
     
Let $\OT$ denote the set of \emph{ordered colored rooted trees}, and
$\OF$ denote the (empty and nonempty) words over the alphabet $\OT$,
henceforth called the set of empty and non-empty \emph{forests}.  It
should be noted that, unlike the classical Butcher theory, the
ordering of the branches in the trees in $\OT$ is important, and
likewise the ordering of the trees within the forest $\OF$.
   
The basic operations involved in building $\OT$ and $\OF$ are:
\begin{itemize}
\item Create the empty forest $\one \in \OF$.
\item Create a longer forest from shorter forests by concatenation,
  $(\omega_1,\omega_2) \mapsto \omega_1\omega_2$.
\item Create a tree from a forest by adding a root node,
  $\signature{B^+}{\OF}{\OT}$.  In the instance where we wish to color
  the nodes using a set of colors $\mathcal{I}$, we introduce an
  indexed family of root adding operations
  $\signature{B^+_i}{\OF}{\OT}$ for all $i\in \mathcal{I}$.  The
  inverse operation whereby we create a forest from a tree by removing
  the root node is written $\signature{B^-}{\OT}{\OF}$. This operation
  extends to $\OF$ by $B^-(\omega_1\omega_2)=
  B^-(\omega_1)B^-(\omega_2)$ and $B^-(\one) = \one$.
\end{itemize}
The total number of forests, with $n$ nodes colored in $i$ different
ways, is defined by modifying the definition of the well known Catalan
numbers
\begin{equation*}
  C_i^n=\frac{i^n}{(n+1)}\binom{2n}{n}, 
  \quad n=0,1,2,\ldots. 
\end{equation*} 
See A000108 in \cite{sloane03tol} for various combinatorial
representations of the Catalan numbers.  For a forest $\omega\in\OF$
we define the \emph{degree}, $\#(\omega)$, as the number of trees in
$\omega$ as:
\begin{equation*}
  \label{eq:degree}
  \begin{split}
  \#(\one) & = 0,\\
  \#(B^+(\omega)) & = 1,\\
  \#(\omega_1\omega_2) & = \#(\omega_1) + \#(\omega_2),
  \end{split}
\end{equation*}
and the \emph{order}, $|\omega|$, as the total number of nodes in all
the trees of $\omega$ as:
\begin{equation*}
  \label{eq:order}
  \begin{split}
  |\one| & = 0,\\
  |B^+(\omega)| & = 1+|\omega|,\\
  |\omega_1\omega_2| & = |\omega_1|+|\omega_2|.
  \end{split}
\end{equation*}
We let $\RN$ denote the linear space of all finite $\RR$-linear
combinations of elements in $\OF$. This vector space is naturally
equipped with an inner product such that all forests are orthogonal,
\begin{equation*}
  \label{eq:orthoword}
  \langle\omega_1,\omega_2\rangle = \left\{
    \begin{array}{cl}
      1, & \mbox{if $\omega_1=\omega_2$},\\
      0, & \mbox{else},
    \end{array}\right.\qquad 
  \mbox{for all $\omega_1,\omega_2\in\OF$}.
\end{equation*}
For $a\in N$ and $\omega\in\OF$, we let $a(\omega)\in\RR$ denote the
coefficient of the forest $\omega$, thus $a$ can be written as a sum
\begin{equation*}
  \label{eq:aseries}
  a = \sum_{\omega\in\OF}a(\omega)\omega,
\end{equation*}
where all but a finite number of terms are zero. The space of all
infinite sums of this kind is denoted $\Nd$ and is the dual space of
$N$, that is
\begin{equation*}
  \Nd = \{\signature{\alpha}{N}{\RR}\colon\,\alpha\,\mbox{linear}\}.
\end{equation*}
We again let $\alpha(a)$ denote the value of $\alpha\in\Nd$ on $a\in
N$.  

The operations $B^+_i$, $B^-$ and concatenation extend to $N$ by
linearity and the distributive law of concatenation, that is
\begin{align*}
  B^+_i(\omega_1+\omega_2) & = B^+_i(\omega_1)+B^+_i(\omega_2),\\
  \omega(\omega_1+\omega_2) & = \omega\omega_1 + \omega\omega_2.
\end{align*}
The vector space $N$ with the concatenation product and the grading
$\#$ forms a  graded associative algebra
$N=\bigoplus_{j\in\ZZ}N_j$, where $N_j$ denotes the linear combination
of forests with $j$ trees.  Alternatively, it is possible to grade
this algebra using $|\cdot |$; in this case, $N_j$ denotes the linear
combination of forests with the same number of nodes.  

Now we introduce a \emph{left grafting product} which has the
algebraic structure of a derivation.
\begin{definition} 
  \label{defn:grafting}
  For $\widetilde{a},a\in N$, define the left grafting $\widetilde{a}[a]\in N$
  by the following recursion formulae, where $\tau\in\OT$ and
  $\omega,\widetilde{\omega}\in\OF$:
  \begin{subequations}
    \label{eq:grafting}
    \begin{equation}
      \label{eq:grafting1}
      \begin{split}
        \tau[\one] & = 0, \\
        \tau[\omega\widetilde{\omega}] & = (\tau[\omega])\widetilde{\omega}+
        \omega(\tau[\widetilde{\omega}]),\\
        \tau[B^+_i(\omega)] & = B^+_i(\tau[\omega])+B^+_i(\tau\omega),
      \end{split}
    \end{equation}
    \begin{equation}
      \label{eq:grafting2}
      \begin{split}
        \one[a] & = a,\\
        (\tau\omega)[a] & = \tau[\omega[a]] - (\tau[\omega])[a].
      \end{split}
    \end{equation}
  \end{subequations}
  The definition of left grafting is extended to the general case
  $\widetilde{a}[a]$ by bilinearity.
\end{definition}
It is useful to understand left grafting directly rather than via the
recursive definition. From~\eqref{eq:grafting1} we verify
that if $\tau\in\OT$ and $\omega\in\OF$ then $\tau[\omega]$ is a
sum of $|\omega|$ words, each word obtained by attaching the root of
$\tau$ with an edge to the left side of a node of $\omega$.
\begin{equation*}
  \aABb\left[\ab\aabABb\aABb\right]=\aaABbb\aabABb\aABb+
  \ab\aaABbabABb\aABb+\ab\aabAaABbBb\aABb+\ab\aaaABbbABb\aABb+
  \ab\aabABb\aaABbabb+\ab\aabABb\aAaABbBb
\end{equation*}
From~\eqref{eq:grafting2} we see that if $\tau_1,\tau_2\in\OT$ and
$\omega\in\OF$ then $(\tau_1\tau_2)[\omega]$ is obtained by first
left-grafting $\tau_2$ to all nodes of $\omega$ and then left-grafting
$\tau_1$ to all the nodes of the resulting expression, except to the
nodes coming originally from $\tau_2$.
\begin{equation*}
  \ab\AB\left[\ab\aABb\right]=\aabABb\aABb+
  \aABb\aabABb+\aABb\aAabBb+\aabb\aABABb+\ab\aabABABb+
  \ab\aABAabBb+\aabb\aAABBb+\ab\aabAABBb+\ab\aAabABBb
\end{equation*}
\begin{lemma}
  \label{lem:1}
  If $\tau_1,\ldots,\tau_k\in\OT$ and $\omega\in\OF$ then
  $(\tau_1\cdots \tau_k)[\omega]$ is a sum of $|\omega|^k$ words
  obtained by, in the order $j=k,k-1,\ldots,1$, attaching the root of
  the tree $\tau_j$ with an edge to the left side of any node in
  $\omega$. In particular we have
  \begin{equation}
    \label{eq:addroot}
    \omega[B^+_i(\one)] = B^+_i(\omega), 
    \qquad\mbox{for all $\omega \in\OF $.}
  \end{equation}
\end{lemma}
Equations~\eqref{eq:grafting1} and~\eqref{eq:grafting2} imply that for
any $d\in N_1$ and $a,b\in N$, we have the Leibniz rule and a
composition rule of the form
\begin{align}
  \label{eq:leibniz}
  d[ab] & = d[a]{b} + a d[{b}] \\
  \label{eq:composit}
  d[a[{b}]] & = da[{b}] + d[a][{b}].
\end{align}
Thus $d$ acts as a first degree {derivation} on $N$.
\begin{definition}
  \label{defn:GL}
  The Grossman--Larson (GL) product $\signature{\opr}{N\tpr N}{N}$ is
  defined as
  \begin{equation*}
    \label{eq:gldef}
    B^+_i\left(\omega\opr\widetilde{\omega}\right) =
    \omega[B^+_i(\widetilde{\omega})],
    \qquad\mbox{for all $\omega,\widetilde{\omega}\in\OF$,}
  \end{equation*}
  and is extended to the general $a\opr\widetilde{a}$ for $a,\widetilde{a}\in
  N$ by linearity.
\end{definition}
Since
$\omega\opr\widetilde{\omega}=B^-\left(\omega[B^+_i(\widetilde{\omega})]\right)$,
the GL product can be understood by adding an invisible root to
$\widetilde{\omega}$ (turning it into a tree), and left-grafting $\omega$
onto all nodes of $B^+_i(\widetilde{\omega})$, including the invisible
root. The root is then removed from each of the resulting trees, with the
GL product resulting in a total of $(|\widetilde{\omega}|+1)^{\#(\omega)}$
forests. Some examples of the grafting product from
Definition~\ref{defn:grafting} and the GL product from
Definition~\ref{defn:GL} are given in Table~\ref{tab:GLproduct}.  
In fact, the theory of Grossman and Larson~\cite{grossman89has} is
formulated on trees, not on forests of trees. To a forest $\omega$ in
our terminology, they add a (proper) root to turn it into a tree. The
definition of the GL product in~\cite{grossman89has} is modified
accordingly.

The GL product is an associative $\#$-graded product,
$\signature{\opr}{N_j\tpr N_k}{N_{j+k}}$, satisfying
for all $a,b,c\in N$, $a\opr\one=\one\opr a$, and 
\begin{align*}
  (a\opr b)[{c}] 
  & = B^-(a[B^+_i(b)])[c]\\
  & = B^-(B^+_i(a[b])+B^+_i(ab))[c]\\
  & = (a[b]+(ab))[c]\\
  & = a[b[{c}]]. 
\end{align*}

\clearpage
\vspace{-5mm}
\begin{table}[t]
  \centering
  \begin{equation*}
    \begin{array}{c@{\,\,}|@{\quad}l|@{\quad}l}
      \hline & & \\
      \omega_1\tpr \omega_2 & \omega[\tau] & \omega\opr\tau \\[2mm]
      \hline & & \\[-1mm]
      \ab\tpr\ab & \aabb &  \ab\ab+\aabb \\[2mm]     
      \ab\ab\tpr\ab & \aababb & \ab\ab\ab+2\ab\aabb+\aababb \\[2mm]
      \ab\tpr\ab\ab & \aabb\ab+\ab\aabb &
      \ab\ab\ab+\aabb\ab+\ab\aabb \\[2mm]
      \aabb\tpr\ab & \aaabbb & \aabb\ab+\aaabbb \\[2mm]
      \ab\tpr\aabb & \aababb+\aaabbb & \ab\aabb+\aababb+\aaabbb \\[2mm]
      \ab\tpr\ab\ab\ab & \aabb\ab\ab+\ab\aabb\ab+\ab\ab\aabb &
      \ab\ab\ab\ab+\aabb\ab\ab+\ab\aabb\ab+\ab\ab\aabb \\[2mm]
      \ab\ab\ab\tpr\ab & \aabababb &
      \ab\ab\ab\ab+3\ab\ab\aabb+3\ab\aababb+\aabababb \\[2mm]
      \ab\aabb\tpr\ab & \aabaabbb &
      \ab\aabb\ab+\ab\aaabbb+\aabb\aabb+\aabaabbb \\[2mm] 
      \aabb\ab\tpr\ab & \aaabbabb &
      \aabb\ab\ab+\aabb\aabb+\ab\aaabbb+\aaabbabb \\[2mm] 
      \ab\tpr\ab\aabb & \aabb\aabb+\ab\aababb+\ab\aaabbb &
      \ab\ab\aabb+\aabb\aabb+\ab\aababb+\ab\aaabbb \\[2mm] 
      \ab\tpr\aabb\ab & \aababb\ab+\aaabbb\ab+\aabb\aabb &
      \ab\aabb\ab+\aababb\ab+\aaabbb\ab+\aabb\aabb \\[2mm] 
      \ab\tpr\aaabbb & \aabaabbb+\aaababbb+\aaaabbbb &
      \ab\aaabbb+\aabaabbb+\aaababbb+\aaaabbbb \\[2mm] 
      \aaabbb\tpr\ab & \aaaabbbb & \aaabbb\ab+\aaaabbbb \\[2mm]
      \ab\tpr\aababb & \aabababb+\aaabbabb+\aabaabbb &
      \ab\aababb+\aabababb+\aaabbabb+\aabaabbb \\[2mm] 
      \aababb\tpr\ab & \aaababbb & \aababb\ab+\aaababbb \\[2mm]
      \ab\ab\tpr\ab\ab & 2\aabb\aabb+\aababb\ab+\ab\aababb &
      \ab\ab\ab\ab+2\ab\aabb\ab+2\ab\ab\aabb+2\aabb\aabb+
      \aababb\ab+\ab\aababb \\[2mm]
      \ab\ab\tpr\aabb & \aabababb+2\aabaabbb+\aaababbb & 
      \ab\ab\aabb+2\ab\aababb+2\ab\aaabbb+\aabababb+2\aabaabbb+\aaababbb
      \\[2mm] 
      \aabb\tpr\ab\ab & \aaabbb\ab+\ab\aaabbb &
      \aabb\ab\ab+\aaabbb\ab+\ab\aaabbb \\[2mm] 
      \aabb\tpr\aabb & \aaabbabb+\aaaabbbb &
      \aabb\aabb+\aaabbabb+\aaaabbbb \\[2mm] \hline  
    \end{array}
  \end{equation*}
  \caption{The left grafting and Grossman--Larson products for all
    forests up to order four. The Grossman--Larson product is the dual of
    the coproduct in $\Hn$ described in
    Section~\ref{sect:hopf}. \label{tab:GLproduct}} 
\end{table}

\clearpage

% ===========================================================================
\subsection{$N$ as a universal object}
\label{sect:universal}
% ===========================================================================

\begin{definition}
  Let $D$ be an associative $\ZZ$-graded algebra
  $D=\bigoplus_{j=0}^\infty D_j$ with associative product $a,b\mapsto
  ab$, a unit $\one$ and grading $\#(D_j)=j$ such that
  $\#(D_jD_k)=j+k$.  We call $D$ a \emph{D-algebra} if it is also
  equipped with a linear derivation $\signature{(\cdot)[\cdot]}{D\tpr
    D}{D}$ such that~\eqref{eq:leibniz} and~\eqref{eq:composit} hold
  for any $d\in D_1$ and any $a,{b}\in D$.
\end{definition}

Define a D-algebra homomorphism as a linear degree preserving map $\F$
between D-algebras satisfying for any $a,b\in D$ :
\begin{align}
  \label{eq:dahom1}  
  \F(ab) & = \F(a)\F(b)\\
  \label{eq:dahom2} 
  \F(a[b]) & = \F(a)[\F(b)].
\end{align}

\begin{proposition}
  Let $N$ be the algebra of forests colored with a set $\mathcal{I}$.
  For any D-algebra $D$ and any map $\signature{i\mapsto
    f_i}{\mathcal{I}}{D_1\subset D}$, there exists a unique
  homomorphism $\signature{\F}{N}{D}$ such that $\F(B^+_i(\one)) =
  f_i$.
\end{proposition}

\begin{proof}
  From~\eqref{eq:addroot} and~\eqref{eq:dahom2} we find
  $\F(B^+_i(\omega)) = \F(\omega)[f_i]$, for any $\omega\in\OF$.
  Together with~\eqref{eq:dahom1} and linearity, this shows that by
  recursion, we can extend $\F$ to a uniquely defined homomorphism
  defined on all of $N$. {\hfill$\Box$ \vskip10pt}
\end{proof}

This shows that $N$ is a universal object, free over the set $\mathcal{I}$, in
the category of $D$-algebras.

% ===========================================================================
\subsection{The algebra of $\Ge$ sections on a manifold}
\label{sect:Gsection}
% ===========================================================================

As an example of a D-algebra, we consider an algebra related to the
numerical Lie group integrators.  Let $\g$ be a Lie algebra of vector
fields on a manifold $\M$ and let
$\signature{\exp}{\g}{\operatorname{Diff}(\M)}$, denote the flow
operator. A basic assumption of numerical Lie group
integrators~\cite{iserles00lgm,munthe-kaas99hor,owren99rkm} is the
existence of a $\g$ which is transitive (i.e.\ spans all tangent
directions in any point on $\mathcal{M}$), and for which the
exponential map can be computed efficiently and exactly. Transitivity
implies that \emph{any} vector field can be written in terms of a
function $\signature{f}{\M}{\g}$. The goal of numerical Lie group
integrators is to approximate the flow of a general differential
equation
\begin{equation}
  \label{eq:ode}
  y'(t) = f(y)(y),
  \qquad\mbox{where $\signature{f}{\M}{\g}$}, 
\end{equation}
by composing exponentials of elements in $\g$. The study of order
conditions for Lie group integrators leads to a need for understanding
the algebraic structure of non-commuting vector fields on $\M$,
generated from $f$.

Elements $V\in\g$ are often called invariant or `frozen' vector fields
on $\M$. These define first-degree invariant differential operators
through the Lie derivative.  Let $\mathcal{V}$ be any normed vector
space and denote $(\M\rightarrow \mathcal{V})$ the set of all smooth
functions from $\M$ to $\mathcal{V}$, called the space of
\emph{$\mathcal{V}$-sections}. For $V\in \g$ and $\psi\in
\left(\M\rightarrow \mathcal{V}\right)$, the Lie derivative,
$V[\psi]\in \left(\M\rightarrow \mathcal{V}\right)$, is defined as
\begin{equation*}
  \label{eq:liederiv}
  V[\psi](p) = \ddto\psi\left(\exp(tV)(p)\right),
  \qquad\mbox{for any point $p\in \M$}.
\end{equation*}
For two elements $V,W\in \g$ we iterate this definition and define the
concatenation $VW$ as the second degree invariant differential
operator $VW[\psi] = V[W[\psi]]$. The linear space spanned by the
$0$-degree identity operator $\one[\psi] = \psi$ and all higher degree
invariant derivations is called the \emph{universal enveloping
  algebra} of $\g$, denoted $\Ge$.  This is a graded algebra with the
concatenation product and degree $\#(\one)=0$, $\#(\g)=1$ and $\#(VW)
= \#(V)+\#(W)$.
  
Given a norm on the vector space $\Ge$, we consider the space of $\Ge$
sections\footnote{Thus $\Ge$ is a trivial vector bundle over $\M$,
  $\g$ a trivial sub-bundle and the tangent bundle $T\M$ is a
  non-trivial sub-bundle of $\g$.} 
$\left(\M\rightarrow\Ge\right)$.  For two sections
$f,g\in\left(\M\rightarrow\Ge\right)$ we define
$f[g]\in\left(\M\rightarrow\Ge\right)$ pointwise from the Lie
derivative as
\begin{equation*}
  \label{eq:sectionderiv}
  f[g](p) = \left(f(p)[g]\right)(p), \qquad p\in\M .
\end{equation*}
Similarly, the concatenation on $\Ge$ is extended pointwise to a
concatenation $fg\in \left(\M\rightarrow\Ge\right)$ as
\begin{equation*}
  (fg)(p) = f(p)g(p), \qquad p\in\M .
\end{equation*}
From these definitions we find:
\begin{lemma}
  \label{lem:dalgebra}
  Let $f\in(\M\rightarrow \g)$ and $g,h\in (\M\rightarrow
  \Ge)$. Then:
  \begin{align*}
    f[gh] & = f[g]h + gf[h], \\
    (f\opr g)[h] \equiv f[g[h]] & = fg[h] + f[g][h].
  \end{align*}
\end{lemma}
\begin{proof} 
  For $p\in \M$ let $V=f(p)\in \g$. Then
  \begin{align*}
    f[gh](p)
    & = \ddto (gh)(\exp(tV)(p))\\
    & = \ddto g(\exp(tV)(p))h(\exp(tV)(p))\\
    & = \ddto g(\exp(tV)(p))h(p) + g(p) h(\exp(tV)(p))\\
    &  = \left(f[g]h + gf[h]\right)(p),\\[2mm]
    (f[g[h]])(p)
    & = \ddto(g(\exp(tf)(p))[h](\exp(tf)(p))\\
    & = \ddto g(p)[h](\exp(tf)(p))+g(\exp(tf)(p))[h](p)\\
    & = (fg)[h](p) + (f[g])[h](p).
  \end{align*}
  {\hfill$\Box$ \vskip10pt}
\end{proof}

Note the difference between $fg$ and $f\opr g$. In the concatenation
the value of $g$ is frozen to ${g}(p)$ before the differentiation with
$f$ is done, whereas in the latter case the spatial variation of $g$
is seen by the differentiation using $f$.

Lemma~\ref{lem:dalgebra} shows that $\left(\M\rightarrow \Ge\right)$
is a $D$-algebra. Thus if we, for every $i\in \mathcal{I}$, pick a
vector field $f_i \in (\M\rightarrow \g)$ then there exists a unique
homomorphism $\signature{\F}{N}{(\M\rightarrow \Ge)}$ such that
$\F(B^+_i(\one)) = f_i$. The images of the trees $\F(\tau)$, for
$\tau\in\OT$, are called the \emph{elementary differentials} in
Butcher's theory (see \cite{butcher63cft}) and the images of the
forests $\F(\omega)$, for $\omega\in\OF$, are called \emph{elementary
  differential operators} in Merson's theory (see \cite{merson57aom}).

% ===========================================================================
\subsection{Elements of Lie--Butcher theory}                  
\label{sect:lbtheory} 
% ===========================================================================

To motivate the algebraic structures of the next section, we briefly
introduce some elements of Lie--Butcher theory. This theory is the
non-commutative generalization of the classical Butcher theory and is
the general foundation behind the construction of order conditions for
Lie group integrators. Various aspects of this theory have been
developed
in~\cite{berland05aso,munthe-kaas95lbt,munthe-kaas98rkm,munthe-kaas03oep,owren99rkm}.
A comprehensive treatment is given in~\cite{munthe-kaas05lbt}.

With the vector space $N$ being defined as the finite $\RR$-linear
combinations of $\OF$, we now let $\Nd$ denote the space of infinite
$\RR$-linear combinations (sequences), or the algebraic dual space of
$N$.  All the operations of the previous paragraphs extend to $\Nd$ by
local finiteness, see~\cite{reutenauer93fla}. For example, the
inner-product in~\eqref{eq:orthoword} is extended to a dual pairing
$\signature{\langle\cdot, \cdot \rangle}{\Nd\xpr N}{\RR}$, where the
computation
$\langle\alpha,b\rangle=\sum_{\omega\in\OF}\alpha(\omega)b(\omega)$ is
always finite, since $b$ is required to be finite.

Consider the homomorphism $\F$ introduced in
Section~\ref{sect:Gsection} extended to a homomorphism of infinite
series $\F \colon \Nd\rightarrow {(\M\rightarrow \Ge)}$, where
$(\M\rightarrow \Ge)$ should now be understood as a space of formal
series. The series might not converge, but all definitions make sense
termwise, and any finite truncation yields a proper $\Ge$-section.  In
classical (commutative) Butcher theory the image of $\alpha\in\Nd$ is
called an S-series, see Murua \cite{murua99fsaI}.  Similarly, we
define an LS-series as an infinite formal series in
$(\M\rightarrow\Ge)$, given by \footnote{See comments at the end of
  Section~\ref{sect:homo} on the chosen normalization.}
\begin{equation}
  \label{eq:lbseries}
  \LS(\alpha) = \sum_{\omega\in\OF}h^{|\omega|}\alpha(\omega)\F(\omega).
\end{equation}
Classical Lie series on manifolds is a generalization of Taylor
series, where the fundamental result is the following `pull-back
formula': Let $f\in (\M\rightarrow\g)$, be a vector field and
$\signature{\exp(f)}{\M}{\M}$, be the $t=1$-flow.  For any $g\in
(\M\rightarrow\g)$ we have (see~\cite{abraham88mta}) that
\begin{equation*}
  g(\exp(f)(p)) = \sum_{j=0}^\infty\frac{1}{j!}f^j[g](p)\equiv
  \operatorname{Exp}(f)[g](p),
\end{equation*}
where $f^0=\one$ and 
$f^j[g] = f[\cdots f[f[g]]\cdots] = (f\opr\cdots\opr f)[g]$. 
Note that if $f = \F(\ab)$, then $f^j = \F(\ab\opr\cdots\opr\ab)$, 
thus the operator exponential
$\operatorname{Exp}(f) = \sum_{j=0}^\infty f^j/j!$ is a LS-series.

Two special cases of LS series are of particular importance: A
LS-series $\LS(\alpha)$ is called \emph{logarithmic} or
\emph{algebra-like} if $\LS(\alpha)\in (\M\rightarrow\g)$ represents a
vector field, and it is is called \emph{exponential} or
\emph{group-like} if $\LS(\alpha)$ is the (formal) operator
exponential of a logarithmic series.  A logarithmic LS-series is the
non-commutative generalization of a B-series, named the
\emph{Lie--Butcher series}~\cite{munthe-kaas05lbt}.

Note that if $\tau_1,\tau_2\in\OT$, then $a = \tau_1\tau_2 -
\tau_2\tau_1$ is a logarithmic series since it represents the
commutator of two vector fields.  More generally, a series
$\alpha\in\Nd$ is logarithmic, if and only if, all its finite components
belong to the free Lie algebra generated by $\OT$. A Hall basis for
this space is characterized in~\cite{munthe-kaas03oep}.
Reutenauer~\cite{reutenauer93fla} presents several alternative
characterizations of logarithmic and exponential series.

We find the characterization in terms of \emph{shuffle products}
particularly useful.  The shuffle product $\signature{\sqcup}{N\tpr
  N}{N}$ is defined for two forests as the summation over all
permutations of the trees in the forests while preserving the ordering
of the trees in each of the initial forests, and is extended to $N\tpr
N$ by linearity.  It can also be recursively defined in the asymmetric
way $\one\sqcup\omega=\omega\sqcup\one=\omega$ for any forest
$\omega\in\OF$, and if $\omega_1=\tau_1\upsilon_1$ and
$\omega_2=\tau_2\upsilon_2$ for $\tau_1,\tau_2\in\OT$ and
$\upsilon_1,\upsilon_2\in\OF$, then
\begin{equation*}
  (\tau_1\upsilon_1)\sqcup(\tau_2\upsilon_2)
  =\tau_1(\upsilon_1\sqcup(\tau_2\upsilon_2))+
  \tau_2((\tau_1\upsilon_1)\sqcup \upsilon_2).
\end{equation*}
The shuffle product is associative and commutative, for all
$\omega_1,\omega_2,\omega_3\in\OF$ we have
\begin{align*}
  (\omega_1\sqcup\omega_2)\sqcup\omega_3 & =
  \omega_1\sqcup(\omega_2\sqcup\omega_3), \\
  \omega_1\sqcup\omega_2 & = \omega_2\sqcup\omega_1.
\end{align*}
Table~\ref{tab:product} gives some simple, but nontrivial, examples of
the shuffle product.
\begin{table}
  \centering
  \begin{equation*}
    \begin{array}{c@{\quad}|@{\quad}l}
      \hline & \\
      \omega_1\tpr \omega_2 & \mu_N(\omega_1\tpr\omega_2) \\[2mm]
      \hline & \\[-1mm]
      \ab\tpr\ab & 2\ab\ab \\[2mm]     
      \ab\tpr\ab\ab & 3\ab\ab\ab  \\[2mm]
      \ab\tpr\aabb & \aabb\ab+\ab\aabb \\[2mm]
      \ab\tpr\ab\ab\ab & 4\ab\ab\ab\ab \\[2mm]
      \ab\tpr\ab\aabb & 2\ab\ab\aabb+\ab\aabb\ab \\[2mm]
      \ab\tpr\aabb\ab & 2\aabb\ab\ab+\ab\aabb\ab \\[2mm]
      \ab\tpr\aaabbb & \ab\aaabbb+\aaabbb\ab \\[2mm]
      \ab\tpr\aababb & \ab\aababb+\aababb\ab \\[2mm]
      \ab\ab\tpr\ab\ab & 6\ab\ab\ab\ab \\[2mm]
      \ab\ab\tpr\aabb & \ab\ab\aabb+\ab\aabb\ab+\aabb\ab\ab \\[2mm]
      \aabb\tpr\aabb & 2\aabb\aabb \\[2mm] \hline
    \end{array}
  \end{equation*}
  \caption{All nontrivial examples of the shuffle product for all
    forests up to and including order four. The shuffle product is
    the dual of the coproduct in the Grossman--Larson Hopf
    algebra. \label{tab:product}} 
\end{table}

\begin{lemma}
  \cite{reutenauer93fla} 
  \label{lem:log_exp}
  A series $\alpha\in\Nd$ is logarithmic if and only if
  \begin{align*}
    \alpha(\one) & = 0\\
    \alpha(\omega_1\sqcup\omega_2) & = 0 
    \qquad\mbox{for all $\omega_1,\omega_2\in \OF\backslash\{\one\}$}.
  \end{align*}
  A series $\alpha\in\Nd$ is exponential if and only if
  \begin{align*}
    \alpha(\one) & = 1 \\
    \alpha(\omega_1\sqcup\omega_2) & = \alpha(\omega_1)\alpha(\omega_2)
    \qquad\mbox{for all $\omega_1,\omega_2\in\OF$}.
  \end{align*}
\end{lemma}

The LS-series of an exponential series $\alpha\in\Nd$ represents
pull-backs, or finite motions on $\M$. They form a group under
composition with the Grossman--Larson product, which is the
generalization of the \emph{Butcher group} to the case of
non-commutative group actions.  Consider the GL product as a linear
operator $\signature{\opr}{\Nd\tpr\Nd}{\Nd}$ defined by
$\opr(\omega_1\tpr\omega_2) = \omega_1\opr\omega_2$. To compute the
composition $\alpha\opr\beta$ for general $\alpha,\beta\in\Nd$, it is
useful to introduce the dual of $\opr$, the linear map
$\signature{\Delta_N}{N}{N\tpr N}$ to be defined in
Section~\ref{sect:recur}. Using Corollary~\ref{cor:dual}, we find
\begin{equation*}
  (\alpha\opr\beta)(\omega) = \langle\opr(\alpha\tpr\beta),\omega\rangle
  = \langle\alpha\tpr\beta,\Delta_N(\omega)\rangle
  = \sum_{\omega_1\tpr\omega_2\in \Delta_N(\omega)}\alpha(\omega_1)\beta(\omega_2) .
\end{equation*}
As an illustrative example, we read from Table~\ref{tab:coproduct}
that
\begin{equation*}
  (\alpha\opr\beta)(\ab\aabb) = \alpha(\ab\aabb)\beta(\one)+
  2\alpha(\ab\ab)\beta(\ab)+\alpha(\ab)\beta(\aabb)+\alpha(\ab)\beta(\ab\ab)
  +\alpha(\one)\beta(\ab\aabb).
\end{equation*}
The inverse in the group is found from the \emph{antipode}, a linear
map $\signature{S_N}{N}{N}$ defined in Section~\ref{sect:hopf}. It can
be shown that the dual of the antipode $\signature{S_N^*}{\Nd}{\Nd}$
also defines the inverse in the group:
\begin{equation*}
  \alpha^{-1}(\omega)=\alpha(S_N(\omega))=\langle\alpha,S_N(\omega)\rangle=\langle 
  S_N^*(\alpha),\omega\rangle,
\end{equation*}
for all exponential $\alpha\in\Nd$. 

% ===========================================================================
\section{Hopf algebras}                  
\label{sect:hopf}
% ===========================================================================

In this section we will study a commutative graded Hopf algebra $\Hn$
of ordered trees. The coproduct $\Delta_N$ in the algebra is defined
by recursion formulae, and at a later stage we will show that this
$\Delta_N$ is the dual of the GL product thereby establishing the
connection between this Hopf algebra and the algebra of the Butcher
group.

We begin by briefly reviewing the definition of a Hopf algebra,
see~\cite{sweedler69ha} for details.  A real associative
\emph{algebra} $\A$ is a real vector space with an associative product
$\signature{\mu}{\A\tpr\A}{\A}$ and a unit $\signature{u}{\RR}{\A}$
such that $\mu(a\tpr u(1)) = \mu(u(1)\tpr a)=a$ for all $a\in\A$.  The
dual of an algebra is called a \emph{coalgebra}, $\C$, which is a
vector space equipped with a coassociative coproduct
$\signature{\Delta}{\C}{\C\tpr \C}$ and counit
$\signature{e}{\C}{\R}$. A \emph{bialgebra} $\B$ is a linear space
which is both an algebra and also a coalgebra structure such that the
coproduct and the counit are compatible with the product, in the sense
that
\begin{align}
  \label{eq:compat1}
  e(\mu(\omega_1\tpr\omega_2)) 
  &= \mu(e(\omega_1)\tpr e(\omega_2)), \\
  \label{eq:compat2}
  \Delta(\mu(\omega_1\tpr\omega_2)) 
  &= (\mu\tpr\mu)(I\tpr T\tpr I) 
  (\Delta(\omega_1)\tpr\Delta(\omega_2)),
\end{align}
where $T(\omega_1\tpr \omega_2) = \omega_2\tpr \omega_1$ is the twist
map.  Let $\End(\B)$ denote all linear maps from $\B$ to itself. We
define the \emph{convolution} $\signature{\star}{\End(\B)}{\End(\B)}$
as
\begin{equation}
  \label{eq:convolution}
  (A\star B)(a)=\mu((A\tpr B)\Delta(a)),
  \qquad\mbox{for $A,B\in\End(\B)$ and $a\in\B$}.
\end{equation}
Let $I\in \End(\B)$ denote the identity matrix.  An \emph{antipode} is
a linear map $S\in\End(\B)$, which is the two-sided inverse of the
identity matrix under convolution, with the antipode satisfying
\begin{equation}
  \label{eq:antipode}
  (I\star S)(a)=(S\star I)(a)=u(e(a)),\qquad\mbox{for all $a\in\B$}.
\end{equation}
\begin{definition}
  A Hopf algebra $\cal{H}$ is a bialgebra equipped with an antipode.
\end{definition}

% ===========================================================================
\subsection{The Hopf algebra of ordered trees}                  
\label{sect:recur}
% ===========================================================================

We will study a particular Hopf algebra based on the vector space of
ordered forests $\RN$, where the coproduct is defined by the following
recursion.
\begin{definition}
  \label{def:coprod}
  Let $\signature{\Delta_N}{N}{N\tpr N}$ be defined by linearity and
  the recursion
  \begin{equation}
    \label{eq:delta-r}
    \begin{split}
      \Delta_N(\one)
      & = \one\tpr\one, \\
      \Delta_N(\omega\tau)
      & = \omega\tau\tpr\one+\Delta_N(\omega)
      \sqcup\cdot(I\tpr B^+_i)\Delta_N(B^-(\tau)),
    \end{split}
  \end{equation}
  where $\tau = B_i^+(\widetilde{\omega})\in\OT$ and
  $\omega,\widetilde{\omega}\in\OF$. The linear operation\\
  $\signature{\sqcup\cdot}{N\tpr N\tpr{N}\tpr{N}}{{N}\tpr{N}}$ is a
  shuffle on the left and concatenation on the right, satisfying
  \begin{equation*}
    (\omega_1\tpr\tau_1)\sqcup\cdot(\omega_2\tpr\tau_2)=
    (\omega_1\sqcup\omega_1)\tpr(\tau_1\tau_2).
  \end{equation*}
\end{definition}
Note that letting $\omega=\one$ yields the special recursion formula
for a tree $\tau$:
\begin{equation*}
  \Delta_N(\tau)
  =\tau\tpr\one+(I\tpr B^+_i)\Delta_N(B^-(\tau)).
\end{equation*}
       
\begin{theorem} 
  \label{thm:hopf} 
  Let $\Hn$ be the vector space $\RN$ with the operations
  \begin{alignat*}{2}
    \mbox{product} & \colon \mu_N(a\tpr b) = a\sqcup b,
    &&\qquad\mbox{(Shuffle product)}\\
    \mbox{coproduct} & \colon \Delta_N,
    &&\qquad\mbox{(Definition~\ref{def:coprod})}\\
    \mbox{unit}       & \colon  u_N(1) = \one,\\
    \mbox{counit} & \colon e_N(\omega) = \left\{\begin{array}{cl}
        1, & \mbox{if $\omega=\one$},\\
        0, & \mbox{else}.
      \end{array}\right. 
  \end{alignat*}
  Then $\Hn$ is a Hopf algebra with an antipode $S_N$ given by the recursion
  \begin{equation}
    \begin{split}
      \label{eq:antipodrecur}
      \S_N(\one) & = \one,\\
      \S_N(\omega\tau) & = -\mu_N\left((\S_N\tpr I)\left(\Delta_N(\omega)\sqcup\cdot
      (I\tpr B^+_i)\Delta_N(B^-(\tau))\right)\right),
    \end{split}
  \end{equation}
  where $\tau = B_i^+(\widetilde{\omega})\in\OT$ and
  $\omega,\widetilde{\omega}\in\OF$.  In particular
  \begin{equation*}
    \S_N(\tau) = -\mu_N\left((\S_N\tpr I)
      (I\tpr B^+_i)\Delta_N(B^-(\tau))\right).
  \end{equation*}
\end{theorem}

\begin{proof}
  The coassociativity of the coalgebra will be established once we
  have shown that the algebra and coalgebra are compatible.  From the
  fact that $e_N(\omega)=0$ for all $\omega\in\OF\backslash\{\one\}$
  and that the shuffle product of two scalars is just standard
  multiplication, we immediately get~\eqref{eq:compat1}.  To
  show~\eqref{eq:compat2} we find it convenient to introduce the
  linear operation
  $\signature{\spr}{{N}\tpr{N}\tpr{N}\tpr{N}}{{N}\tpr{N}}$ with the
  shuffle product both on the left and the right, satisfying
  \begin{equation*}
    (\upsilon_1\tpr\nu_1)\spr(\upsilon_2\tpr\nu_2)=
    (\upsilon_1\sqcup\upsilon_2)\tpr(\nu_1\sqcup\nu_2).
  \end{equation*}
  The compatibility condition \eqref{eq:compat2} is now equivalent to
  \begin{equation}
    \label{eq:coprod-r}
    \Delta_N(\omega_1\sqcup\omega_2)=\Delta_N(\omega_1)\spr\Delta_N(\omega_2).
  \end{equation}
  To simplify the notation we use the fact that
  \begin{equation*}
    \Delta_N(\omega\tau)=\omega\tau\tpr\one+\overline{\Delta}_N(\omega\tau).
  \end{equation*}
  Let $\omega_1=\tau_1\upsilon_1$ and $\omega_2=\tau_2\upsilon_2$
  for $\tau_1,\tau_2\in\OT$ and $\upsilon_1,\upsilon_2\in\OF$, now
  using the recursive definition of the shuffle product and
  substituting the expression for the coproduct of an ordered forest
  yields
  \begin{align*}
    &\Delta_N(\upsilon_1\tau_1\sqcup \upsilon_2\tau_2)
    =\Delta_N((\upsilon_1\sqcup\tau_2\upsilon_2)\tau_1)+
    \Delta((\tau_1\upsilon_1\sqcup u_2)\tau_2)\\
    &=(\Delta_N(\upsilon_1)\spr\Delta_N(\upsilon_2\tau_2))\sqcup
    \cdot\overline{\Delta}_N(\tau_1)+
    (\Delta_N(\upsilon_1\tau_1)\spr\Delta_N(\upsilon_2))\sqcup
    \cdot\overline{\Delta}_N(\tau_2)\\
    &+(\upsilon_1\tau_1\sqcup \upsilon_2)\tau_2\tpr\one+
    (\upsilon_1\sqcup \upsilon_2\tau_2)\tau_1\tpr\one \\
    &=(\Delta_N(\upsilon_1)\spr\overline{\Delta}_N(\upsilon_2\tau_2))\sqcup
    \cdot\overline{\Delta}_N(\tau_1)+
    (\Delta_N(\upsilon_1\tau_1)\spr\overline{\Delta}_N(\upsilon_2))\sqcup
    \cdot\overline{\Delta}_N(\tau_2)\\
    &+(\Delta_N(\upsilon_1)\spr\upsilon_2\tau_2\tpr\one)\sqcup\cdot
    \overline{\Delta}_N(\tau_1)+(\Delta_N(\upsilon_2)\spr\upsilon_1\tau_1
    \tpr\one)\sqcup\cdot\overline{\Delta}_N(\tau_2)\\
    &+(\upsilon_1\tau_1\sqcup \upsilon_2\tau_2)\tpr\one\\
    &=\upsilon_1\tau_1\sqcup \upsilon_2\tau_2\tpr\one+
    \overline{\Delta}_N(\upsilon_1\tau_1)\spr(\upsilon_2\tau_2\tpr\one)+
    \overline{\Delta}_N(\upsilon_2\tau_2)\spr(\upsilon_1\tau_1\tpr\one)\\
    &+\overline{\Delta}_N(\upsilon_1\tau_1)\spr\overline{\Delta}_N
    (\upsilon_2\tau_2)\\
    &=\Delta_N(\upsilon_1\tau_1)\spr\Delta(\upsilon_2\tau_2).
  \end{align*}
  We now have the tools needed to show that the coalgebra is
  coassociative, which follows from the fact that the coproduct
  $\Delta_N$ satisfies
  \begin{align*}
    (I\tpr\Delta_N)\Delta_N(\omega_1\sqcup\omega_2)
    &=(I\tpr\Delta_N)(\Delta_N(\omega_1)\spr\Delta_N(\omega_2)) \\
    &=(I\tpr\Delta_N)\Delta_N(\omega_1)\spr(I\tpr\Delta_N)\Delta_N(\omega_2) \\
    &=(\Delta_N\tpr I)\Delta_N(\omega_1)\spr(\Delta_N\tpr I)\Delta_N(\omega_2) \\
    &=(\Delta_N\tpr I)\Delta_N(\omega_1\sqcup\omega_2).
  \end{align*}
  Thus we have established the structure of a bialgebra.  Substituting
  Definition~\ref{def:coprod} in~\eqref{eq:antipode} yields the
  recursion for the antipode~\eqref{eq:antipodrecur}.  {\hfill$\Box$
    \vskip10pt}
\end{proof}

Both the definition of the coproduct $\Delta_N$ and thus the antipode
$S_N$ are recursive and difficult to use in practice. To develop
non-recursive formulae for these, it is first necessary to define
certain cutting operations.

\begin{definition}
  For a given forest $\omega\in\OF$, a \emph{parent} is any node $p$
  with at least one branch growing from that node and the
  \emph{children} are the nodes branching from $p$. Let $p_c$ denote
  the number of children of $p$.  Cutting off a child node equates to
  removing the edge connecting the child to its parent.
  \begin{itemize}
  \item A \emph{nodal left cut} of degree $c$ is a cut where the $c$
    leftmost children of a given parent node $p$ are cut off.  We can
    write a nodal left cut as $\ell_p(c)$ where $0<c\leq p_c$.  The
    cut splits a forest $\omega$ into two sub-forests,
    $P^{\ell_p(c)}(\omega)$ and $R^{\ell_p(c)}(\omega)$, where $P$ is
    the part cut off, with the forest containing the $c$ leftmost
    children of $p$ as root nodes and $R$ is the remaining bottom part
    of $\omega$.
  \item A \emph{left cut} is a collection of $0\leq k$ nodal left cuts
    $\ell = \left\{\ell_{p_i}(c_i)\right\}_{i=1}^k$, where
    $\{p_i\}_{i=1}^k$ are distinct nodes of $\omega$.  This splits
    $\omega$ in $k$ cut-off forests $\{\omega_i\}_{i=1}^k$ and a
    remaining forest $R^\ell(\omega)$, where $\omega_i$ is the forest
    containing the nodes connected to the $c_i$ leftmost children of
    $p_i$ and $R^\ell(\omega)$ is the forest of the nodes connected to
    the original root nodes.  We define $P^{\ell}(\omega)\in N$ as
    \begin{equation}
      \label{eq:P}
      P^{\ell}(\omega) = \omega_1\sqcup\omega_2\sqcup\cdots \sqcup \omega_k.
    \end{equation}
    Note that the definition of a left cut includes the case $k=0$,
    called the \emph{empty cut}, where $R^\ell(\omega) = \omega$ and
    we define $P^\ell(\omega) = \one$.
  \item An \emph{admissible left cut} is a left cut, containing the
    restriction that any path from a node in $\omega$ to the
    corresponding root is cut no more than once.
  \end{itemize}
\end{definition}

We denote by LC, NLC and ALC the set of all left cuts, nodal left cuts
and admissible left cuts.  To define the coproduct we need to slightly
extend the definition of an admissible left cut, which we choose to
call a \emph{full admissible left cut}. The full admissible left cuts
of $\omega\in\OF$ are obtained by adding an (invisible) root node to
form the tree $\tau = B_i^+(\omega)\in\OT$, applying an admissible
left cut on $\tau$, and finally removing the invisible root node
again. We denote by FALC the set of all full admissible left
cuts. Thus $\mathrm{FALC}(\omega) = \mathrm{ALC}(\tau)$ and for any
$\ell\in\mathrm{FALC}(\omega)$ we have $P^\ell(\omega) = P^\ell(\tau)$
and $R^\ell(\omega) = B^-(R^\ell(\tau))$.  Note that
$\mathrm{FALC}(\omega)$ contains the `cut everything', where $k=1$,
$P^\ell(\omega) = \omega$ and $R^\ell(\omega) = \one$, as well as the
empty cut with $k=0$, $P^\ell(\omega) = \one$ and $R^\ell(\omega) =
\omega$. It is useful to note that the order in which the cuts are
performed does not affect $P^\ell(\omega)$ or $R^\ell(\omega)$. The
order of the cuts is taken care of by the use of the shuffle product
in the definition of $P^\ell(\omega)$ given by equation \eqref{eq:P}.

As an example, we list all the cuttings of an example tree in
Table~\ref{tab:cutexample}.
\begin{table}[ht]
  \centering
  \begin{equation*}
    \begin{array}{c@{\quad\!}|@{\quad\!}c@{\quad\!}|@{\quad\!}c@{\quad\!}|
        @{\quad\!}c@{\quad\!}|@{\quad\!}c@{\quad\!}||@{\quad\!}c@{\quad\!}|
        @{\quad\!}c@{\quad\!}|@{\quad\!}c@{\quad\!}|@{\quad\!}c@{\quad\!}|
        @{\quad\!}c}
      \hline & & & & & & & & \\
      i & \ell_i & k & P^{\ell_i}(\tau) & R^{\ell_i}(\tau) & 
      i & \ell_i & k & P^{\ell_i}(\tau) & R^{\ell_i}(\tau) \\[2mm]  
      \hline  & & & & & & & & \\[-1mm]
      0 & \trxl & 0 & \one & \trxl & 7 & \trxg & 2 & \aababb\sqcup\ab
      & \ab \\[2mm]
      1 & \trxa & 1 & \ab & \aaababbb & 8 & \trxh & 2 &
      \aaabbb\sqcup\ab & \ab \\[2mm]
      2 & \trxb & 1 & \ab & \aaaabbbb & 9 & \trxi & 2 &
      \ab\sqcup\ab\aabb & \ab \\[2mm]
      3 & \trxc & 1 & \ab\aabb & \aabb & 10 & \trxj & 3 &
      \aabb\sqcup\ab\sqcup\ab & \ab \\[2mm]
      4 & \trxd & 1 & \aabaabbb & \ab & 11 & \trxk & 3 &
      \ab\sqcup\ab\ab\sqcup\ab & \ab \\[2mm]
      5 & \trxe & 2 & \ab\sqcup\ab & \aaabbb & 12 & \trxl & 1 &
      \aaabaabbbb & \one \\[2mm]
      6 & \trxf & 2 & \ab\ab\sqcup\ab & \aabb & & & & & \\[2mm] \hline
    \end{array}
  \end{equation*}
  \caption{The cuts $\ell_i$ of an example tree $\tau$, where
    $\ell_{12}$ is the `cut everything' full cut. Thus
    $\mathrm{NLC}(\tau)=\{\ell_1,\ldots,\ell_4\}$,
    $\mathrm{LC}(\tau)=\{\ell_0,\ldots,\ell_{11}\}$,
    $\mathrm{ALC}(\tau) =\{\ell_0,\ldots,\ell_5\}$ and
    $\mathrm{FALC}(\tau)=\{\ell_0,\ldots,\ell_5\}\cup
    \{\ell_{12}\}$. \label{tab:cutexample}} 
\end{table}

\begin{proposition}
  The coproduct $\Delta_N$ of $\Hn$ is non-recursively
  defined as
  \begin{equation}
    \label{eq:delta-nr}
    \Delta_N(\omega)= 
    \sum_{\ell\in\mathrm{FALC}(\omega)} P^\ell(\omega)\tpr R^\ell(\omega).
  \end{equation}
\end{proposition}
\begin{proof}
  To prove that the recursive definition \eqref{eq:delta-r} and the
  non-recursive definition \eqref{eq:delta-nr} of the coproduct are
  identical, an induction argument on the number of vertices is used.
  First recall that $\mathrm{FALC}(B^-(\tau))=\mathrm{ALC}(\tau)$, and
  for any $\ell\in\mathrm{FALC}(B^-(\tau))$ we have $P^\ell(B^-(\tau))
  = P^\ell(\tau)$ and $R^\ell(B^-(\tau)) = B^-(R^\ell(\tau))$, this
  implies that
  \begin{equation*}
    \sum_{\jmath\in\mathrm{FALC}(B^-(\tau))}
      P^\jmath(B^-(\tau))\tpr B_i^+(R^\jmath(B^-(\tau)))=
    \sum_{\jmath\in\mathrm{ALC}(\tau)}P^\jmath(\tau)\tpr
      R^\jmath(\tau).
  \end{equation*}
  Using this fact, the coproduct now takes the form
  \begin{align*}
    \Delta_N(\omega\tau)
    &=\omega\tau\tpr\one+\left(\sum_{\ell\in\mathrm{FALC}(\omega)}
      P^\ell(\omega)\tpr R^\ell(\omega)\right)\!\sqcup\cdot\!
    \left(\sum_{\jmath\in\mathrm{ALC}(\tau)}P^\jmath(\tau)\tpr
      R^\jmath(\tau)\right)\\
    &=\omega\tau\tpr\one+\sum_{\ell\in\mathrm{FALC}(\omega)}\sum_{\jmath\in\mathrm{ALC}(\tau)}
    P^\ell(\omega)\sqcup P^\jmath(\tau)\tpr R^\ell(\omega)R^\jmath(\tau)\\
    &=\sum_{\ell\in\mathrm{FALC}(\omega\tau)}P^\ell(\omega\tau)\tpr R^\ell(\omega\tau).
  \end{align*} 
  The last equality is true because the sum over
  $\ell\in\mathrm{FALC}(\omega)$ and $\jmath\in\mathrm{ALC}(\tau)$ is
  equivalent to the sum over $\ell\in\mathrm{FALC}(\omega\tau)$ except
  for the `cut everything' cut which is equal to the term
  $\omega\tau\tpr\one$. {\hfill$\Box$ \vskip10pt}
\end{proof}

\begin{corollary}
  \label{cor:dual}
  The dual of the coproduct $\Delta_N$ is $\opr$ the Grossman--Larson
  product, that is for any $\omega\in N$ and $\omega_1,\omega_2\in
  \Nd$ we have
 \begin{equation*}
   \langle\omega_1\opr\omega_2,\omega\rangle =
   \langle\omega_1\tpr\omega_2,\Delta_N(\omega)\rangle.
 \end{equation*}
\end{corollary}
\begin{proof}
  If the sum in \eqref{eq:delta-nr} had been over ALC instead of FALC,
  then the dual would have been the left grafting. To see this, we use
  the characterization of left grafting in Lemma~\ref{lem:1},
  and observe that the nodal left cut corresponds to the dual
  operation of attaching a number of trees in a given order to a
  common node, while the shuffles in $P^\ell(\omega)$ corresponds to
  the dual operation of attaching the forests in all possible ways to
  different nodes.  From Definition~\ref{defn:GL}, we see that when
  the sum is extended from ALC to FALC, then we obtain the dual of the
  Grossman--Larson product.  {\hfill$\Box$ \vskip10pt}
\end{proof}

To present a non-recursive definition of the antipode, we define the
{reversal map} ${S_F}\colon N\rightarrow N$ as
\begin{equation}
  \label{eq:rho}
  \begin{split}
    S_F(\one)&=\one,\\
    S_F(\tau_1\tau_2\cdots \tau_j)&= (-1)^j \tau_j\tau_{j-1}\cdots
    \tau_1, \quad\mbox{for all $\tau_1\cdots\tau_j \in \OF$},
  \end{split}
\end{equation}
extended to $N$ by linearity. Thus ${S_F}$ is the unique
anti-automorphism of the concatenation algebra which sends
$\tau\mapsto -\tau$.

\begin{proposition}
  \label{prop:antipode_nr}
  The antipode $S_N$ of $\Hn$ is non-recursively defined as
  \begin{equation}
    \label{eq:antipode-nr}
    S_N(\omega) = {S_F}\left(\sum_{\ell\in\mathrm{LC}(\omega)}  
      P^\ell(\omega)\sqcup R^\ell(\omega)\right).
  \end{equation}
\end{proposition}

\begin{proof}
  In order to prove this result, we need some results about a Hopf
  algebraic structure of the Free Associative Algebra
  (FAA)~\cite{reutenauer93fla}.  Given an alphabet $A$, FAA is the
  vector space formed by taking all finite linear combinations of
  words over $A$. In our case, the alphabet is $\OT$, the words are
  $\OF$ and the vector space is $N$.  A Hopf algebraic structure $\Hf$
  is obtained by taking the product $\mu_F=\mu_N$ as the shuffle
  product and the coproduct $\Delta_F$ defined as the dual of the
  concatenation product. The antipode is the map $S_F$ defined
  in~\eqref{eq:rho}. We need a characterization of $\Delta_F$
  and $S_F$ in terms of cutting operations. For a $\omega\in \OF$
  let the set of Word Cuts (WC) be a simple cut $\ell$ which splits a
  word $\omega$ into two parts $\omega_1 = P^\ell(\omega)$ and
  $\omega_2 = R^\ell(\omega)$ such that $\omega =
  \omega_1\omega_2$. WC contains both the empty cut where
  $P^\ell(\omega) = \one$, $R^\ell(\omega) = \omega$ and cut
  everything where $P^\ell(\omega) = \omega$, $R^\ell(\omega) =\one$.
  Note that the difference between ALC and FALC is that FALC may
  contain a nonempty cut from WC. A direct definition of $\Delta_F$ is
  \begin{equation*}
    \Delta_F(\omega) = \sum_{\ell\in\mathrm{WC}}
    P^\ell(\omega)\tpr R^\ell(\omega),
    \qquad\mbox{for all $\omega\in \OF$.}
  \end{equation*}
  From~\eqref{eq:antipode} we find for
  $\omega\in\OF\backslash\{\one\}$ that
  \begin{equation*}
    0 = (S_F\star I)(\omega) = \mu_F((S_F\tpr I)\Delta_F(\omega)) = 
    \sum_{\ell\in \mathrm{WC}(\omega)} S_F(P^\ell(\omega)) \sqcup
    R^{\ell}(\omega).
  \end{equation*}
  Thus we find a recursive definition of the antipode $S_F$
  \begin{equation}
    \label{eq:sf_recur} 
    \begin{split}
      S_F(\one) &= \one,\\
      S_F(\omega) &= -\sum_{\ell\in\mathrm{WC}(\omega)\backslash
        \mathrm{c.e.}}S_F(P^\ell(\omega)) \sqcup R^{\ell}(\omega),
    \end{split}
  \end{equation}
  where \emph{c.e.} denotes \emph{cut everything}.  Now we repeat the
  same computation for $S_N$, using~\eqref{eq:delta-nr}. This gives
  the recursive definition of the antipode $S_N$
 \begin{equation}
    \label{eq:sn_recur} 
    \begin{split}
      S_N(\one) &= \one,\\
      S_N(\omega) &= -\sum_{\ell\in\mathrm{FALC}(\omega)\backslash
        \mathrm{c.e.}} S_N(P^\ell(\omega)) \sqcup R^{\ell}(\omega).
    \end{split}
  \end{equation}
  We prove~\eqref{eq:antipode-nr} by induction on the number of
  nodes. Plugging~\eqref{eq:antipode-nr} into~\eqref{eq:sn_recur}, we
  find for $\omega\in\OF\backslash\{\one\}$ that
  \begin{align*}
    S_N(\omega) &= -\sum_{\ell\in\mathrm{FALC}(\omega)\backslash
      \mathrm{c.e.}}{S_F}\left(\sum_{j\in\mathrm{LC}(P^\ell(\omega))}
      P^j(P^\ell(\omega))\sqcup R^j(P^\ell(\omega)) \right)\sqcup R^{\ell}(\omega)\\
    &=-\sum_{j\in\mathrm{LC}(\omega)}\sum_{\ell\in\mathrm{WC}(\omega)\backslash
      \mathrm{c.e.}}{S_F}\left( P^j(\omega)\sqcup P^\ell(R^j(\omega))
    \right)\sqcup R^{\ell}(R^j(\omega)) \\
    &= \sum_{j\in\mathrm{LC}(\omega)}\!S_F\left(P^j(\omega)\right)
    \sqcup\left(-\!\sum_{\ell\in\mathrm{WC}(\omega)\backslash
        \mathrm{c.e.}}\!
      S_F(P^\ell(R^j(\omega)))\sqcup R^{\ell}(R^j(\omega))  \right)\\
    &= \sum_{j\in\mathrm{LC}(\omega)} {S_F}\left(P^j(\omega)\right)
    \sqcup S_F(R^j(\omega)) \\
    &= S_F\left(\sum_{j\in\mathrm{LC}(\omega)}
      P^j(\omega)\sqcup R^j(\omega) \right).
  \end{align*}
  We have used the relation $S_F(\omega_1\sqcup\omega_2) =
  S_F(\omega_1)\sqcup S_F(\omega_2)$ see Corollary \ref{coro:connect}
  and the recursion~\eqref{eq:sf_recur}, as well as a careful
  replacement of the summations over $\ell\in\mathrm{FALC}(\omega)$
  and $\jmath\in\mathrm{LC}(P^\ell(\omega))$ with an equivalent sum
  over $\jmath\in\mathrm{LC}(\omega)$ and
  $\ell\in\mathrm{WC}(\omega)$.  {\hfill$\Box$ \vskip10pt}
\end{proof}

As an example, we compute $\Delta_N$ and $S_N$ for the word
$\omega=\ab\aaabbb$. The cuts LC and FALC are shown in
Table~\ref{tab:directex}.  From the direct formulae we find
$\Delta_N(\omega)$ and $\S_N(\omega)$ as listed in
Table~\ref{tab:coproduct} and Table~\ref{tab:antipode}.
\begin{table}[ht]
  \centering
  \begin{equation*}
    \begin{array}{c@{\quad}|@{\quad}c@{\quad}|@{\quad}c@{\quad}|@{\quad}
        c@{\quad}|@{\quad}c@{\quad}||@{\quad}c@{\quad}|@{\quad}c@{\quad}|
        @{\quad}c@{\quad}|@{\quad}c@{\quad}|@{\quad}c}
      \hline & & & & & & & & & \\
      i & \ell_i & k & P^{\ell_i}(\omega) & R^{\ell_i}(\omega) &  
      i & \ell_i & k & P^{\ell_i}(\omega) & R^{\ell_i}(\omega) \\[2mm]  
      \hline  & & & & & & & & & \\[-1mm]
      0 & \ab\aaabbb & 0 & \one & \ab\aaabbb & 
      4 & \ab\aaabbb  & 1 & \ab & \aaabbb \\[2mm]
      1 & \ab\aaabbbz & 1 & \ab & \ab\aabb & 
      5 & \ab\aaabbbz & 2 & \ab\sqcup\ab & \aabb \\[2mm]
      2 & \ab\aaabbby & 1 & \aabb & \ab\ab &
      6 & \ab\aaabbby & 2 & \ab\sqcup\aabb & \ab \\ [2mm]
      3 & \ab\aaabbbx & 2 & \ab\sqcup\ab & \ab\ab & 
      7 & \ab\aaabbb  & 1 & \ab\aaabbb & \one \\[2mm] \hline
  \end{array}
  \end{equation*}
  \caption{Cuts $\ell_i$ of an example word $\omega$.  
    The cuts $\{\ell_{4},\ell_5,\ell_6\}$ are full cuts where the
    leftmost child of the invisible root is cut and $\ell_7$ is the
    full cut where both the children of the invisible root are cut.
    Thus $\mathrm{NLC}(\omega)=\{\ell_1,\ell_{2}\}$,
    $\mathrm{LC}(\omega)=\{\ell_0,\ldots,\ell_{3}\}$, 
    $\mathrm{ALC}(\omega)=\{\ell_0,\ldots,\ell_{2}\}$, 
    $\mathrm{FALC}(\omega)=\{\ell_0,\ldots,\ell_2\}\cup
    \{\ell_{4},\ldots,\ell_{7}\}$ and $\mathrm{WC}(\omega)= 
    \{\ell_0,\ell_4,\ell_7\}$. \label{tab:directex}}  
\end{table}

We complete this section by listing some well known but very useful
relations of Hopf algebras, see Sweedler \cite{sweedler69ha} for
further details. 
\begin{corollary}
  \label{coro:connect}
  \cite{sweedler69ha}
  Given $\Hn$ a Hopf algebra, with product $\sqcup$, coproduct
  $\Delta_N$ and antipode $S_N$, then for all
  $\omega_1,\omega_2\in\OF$  
  \begin{align*}
    (S_N\tpr S_N)\Delta_N(\omega_1) & = \Delta_N(S_N(\omega_1)),\\
    S_N(\omega_1)\sqcup S_N(\omega_2) & =
    S_N(\omega_1\sqcup\omega_2).
  \end{align*}
  If $\Hn$ is either commutative or cocommutative, then
  $S_N(S_N(\omega))=\omega$ for all $\omega\in\OF$.
 \end{corollary}

% ===========================================================================
\section{Hopf algebras related to $\Hn$}                  
\label{sect:homo}
% ===========================================================================

There are two interesting commutative graded Hopf sub-algebras of
$\Hn$ obtained by restricting from the set of ordered rooted trees to
either the set of tall trees (that is trees where each parent has one
child) or bushy trees (that is trees where there is only one parent).
These Hopf sub-algebras are useful respectively for determining the
order conditions for the problem \eqref{eq:ode}, when $f(y)$ is
constant, or when the numerical scheme has high stage order.  

In this section we will also show that the Hopf algebra $\Hc$ of
Butcher, based on unordered trees, can be identified as a
sub-algebra of $\Hn$.  Finally we find that the Hopf algebra $\Hf$ of
the Free Associative Algebra is related to $\Hn$ through the operation
of \emph{freezing} vector fields, which can be defined as a quotient
construction on $\Hn$.

% ===========================================================================
\subsection{Connections to the Butcher theory}
% ===========================================================================

Let $\UT$ denote all unordered trees and $\UF$ denote all unordered
forests, defined as the set of all empty or non-empty unordered words
over the alphabet $\UT$.  Recall from \cite{connes98har} the following
definition.
\begin{definition}
  Given the real vector space $\RC$, denote the Hopf algebra of
  unordered forests as $\Hc=(C,\mu_C,u_C,\Delta_C,e_C,S_C)$. The
  product $\signature{\mu_C}{C\tpr C}{C}$ is defined as the
  (commutative) concatenation
  \begin{equation*}
    \mu_C(\omega_1\tpr\omega_2)=\omega_1\omega_2.
  \end{equation*}
  The unit element $\signature{u_C}{\R}{C}$, is given by
  $u_C(1)=\one$.  The coproduct $\signature{\Delta_C}{C}{C\tpr C}$ is
  defined by linearity and for any
  $\tau=B_i^+(\widetilde{\omega})\in\UT$ and
  $\omega,\widetilde{\omega}\in\UF$ by the recursion
  \begin{equation}
    \label{eq:DcRecur}
    \begin{split}
      \Delta_C(\one) & = \one\tpr\one,\\
      \Delta_C(\tau) & = \tau\tpr\one+(I\tpr B^+_i)\Delta_C(B^-(\tau)),\\
      \Delta_C(\omega\tau) & = \Delta_C(\omega)\Delta_C(\tau).
    \end{split}
  \end{equation}
  The counit $\signature{e_C}{C}{\R}$ is defined by $e_N(\one)=1$ and
  $e_N(\omega)=0$ for $\omega\in \OF\backslash\{\one\}$.  The
  antipode $\signature{S_C}{C}{C}$ is, as usual, the two-sided inverse
  of the convolution in $\Hc$, see \cite{connes98har} for details.
\end{definition}

The main tool used to provide the relationship between the Hopf
algebras of ordered and unordered forests is the symmetrization
operator defined below.
\begin{definition}
  \label{defn:symm}
  The symmetrization operator $\signature{\Omega}{N}{N}$ is a map
  defined by linearity and the relations
  \begin{align*}
    \Omega(\one) & = \one,\\
    \Omega(\omega\tau) & = \Omega(\omega)\sqcup\Omega(\tau),\\
    \Omega(B^+_i(\omega)) & = B^+_i(\Omega(\omega)).
  \end{align*}
\end{definition}

The shuffle product permutes the trees in a forest in all possible
ways, and the symmetrization of a tree is a recursive splitting in
sums over all permutations of the branches.  The symmetrization
defines an equivalence relation on $\OF$, that is, if
\begin{equation*}
  \label{eq:eqforests}
  \Omega(\omega_1)=\Omega(\omega_2) 
  \quad\Longleftrightarrow\quad \omega_1\sim\omega_2.
\end{equation*}
Thus $\omega_1\sim\omega_2$ if and only if $\omega_2$ can be obtained
from $\omega_1$ by permuting the order of the trees in the forest and
the order of the branches of the trees.  We see that an alternative
characterization of $\Omega$ is
\begin{equation*}
  \label{eq:symm4}
  \Omega(\omega)  = \sigma(\omega)\sum_{\stackrel{\widetilde{\omega}\in\OF}
    {\widetilde{\omega}\sim\omega}}\widetilde{\omega} .
\end{equation*}
The integer $\sigma(\omega)$ is the classical symmetry coefficient,
defined for trees and forests as
\begin{align*}
  \sigma(\one) & = 1, \\
  \sigma(\tau_1\tau_2\cdots\tau_k) & =
  \sigma(\tau_1)\cdots\sigma(\tau_k) \mu_1!\mu_2!\ldots, \\
  \sigma(B_i^+(\tau_1\cdots\tau_k)) &= \sigma(\tau_1\cdots\tau_k),
\end{align*}
where the integers $\mu_1,\mu_2,\ldots$ count the number of equivalent
trees among $\tau_1,\ldots,\tau_k$. In other words, if we consider the
full group of all possible permutations of trees and branches acting
on a forest $\omega\in \OF$, then $\sigma(\omega)$ is the size of the
isotropy subgroup i.e.\ the number of permutations leaving $\omega$
invariant. The total number of permutations acting on a given forest
$\omega$ is given by the integer $\pi(\omega)$ defined as
\begin{align*}
  \pi(\one) & = 1, \\
  \pi(\tau_1\tau_2\cdots\tau_k) & =
  k! \sigma(\tau_1)\cdots\sigma(\tau_k) , \\
  \pi(B_i^+(\tau_1\cdots\tau_k)) &= \pi(\tau_1\cdots\tau_k).
\end{align*}
Note that once a forest $\omega$ is symmetrized, then another
application of the symmetrization yields the scaling
\begin{equation}
  \label{eq:scaling}
  \Omega(\Omega(\omega)) = \pi(\omega)\Omega(\omega).
\end{equation}

Let $\UF$ be the unordered forests. Clearly, there is a 1--1
correspondence between unordered forests and equivalence classes of
ordered forests, thus there is a natural isomorphism $\UF \simeq
\OF/\sim$.  Through this identification, we can interpret $\Omega$ as
an injection $\signature{\Omega}{C}{N}$ where $\RC$ and
$\RN$. From~\eqref{eq:scaling} we see that the map
$\signature{\Omega^{-1}}{N}{C}$ defined as
\begin{equation*}
  \Omega^{-1}(a) = \sum_{\omega\in\OF}\frac{a(\omega)}{\pi(\omega)}
  \operatorname{forget}(\omega), 
\end{equation*}
where $\signature{\operatorname{forget}}{N}{C}$ is the natural
identification of an ordered forest with the corresponding unordered
forest, defines a left-sided inverse $\Omega^{-1}(\Omega(b)) = b$ for
all $b\in C$.

\begin{theorem}
  \label{thm:Hopfhomo}
  The symmetrization operator $\signature{\Omega}{C}{N}$ defines an
  injective Hopf algebra homomorphism from the Hopf algebra $\Hc$ of
  unordered forests into the Hopf algebra $\Hn$ of ordered forests.
\end{theorem}
\begin{proof}
  A Hopf algebra homomorphism is a bialgebra homomorphism, which is a
  linear map that is both an algebra and a coalgebra homomorphism.
  ${\Omega}$ is an algebra homomorphism if
  \begin{align*}
    \Omega(u_C(1)) & = u_N(1),\\
    \mu_N(\Omega(\omega_1)\tpr\Omega(\omega_2)) & =
    \Omega(\mu_C(\omega_1\tpr\omega_2)).
  \end{align*}
  These conditions are automatically satisfied by Definition
  \ref{defn:symm}. ${\Omega}$ is a coalgebra homomorphism if
  \begin{align}
    \label{eq:omegarelation1}
    e_N(\Omega(\omega)) & = e_C(\omega),\\
    \label{eq:omegarelation2}
    \Delta_N(\Omega(\omega)) & = (\Omega\tpr\Omega)\Delta_C(\omega).
  \end{align}
  The first condition \eqref{eq:omegarelation1} follows immediately
  given that the counits are only non-zero when the argument is the
  empty forest.
  
  The second relation follows using an induction argument. First we
  need to establish a useful relationship between $\Omega$ and
  $\Delta_C$. Using the \emph{sumless Sweedler notation}
  \begin{equation*}
    \Delta(\omega)=\sum_{\omega_1\tpr\omega_2\in\Delta(\omega)}
    \omega_1\tpr\omega_2=\sum_i\omega_{(1)}^{(i)}\tpr\omega_{(2)}^{(i)}\equiv
    \sum_i\omega_{(1)}\tpr\omega_{(2)}\equiv\omega_{(1)}\tpr\omega_{(2)},
  \end{equation*}
  we find:
  \begin{align}
    \label{eq:relate}
    (\Omega\tpr\Omega)\Delta_C(\omega)\spr(\Omega\tpr\Omega)
    \Delta_C(\tau)
    &=(\Omega\tpr\Omega)\left(\omega_{(1)}\tpr\omega_{(2)}\right)
    \spr(\Omega\tpr\Omega)\left(\tau_{(1)}\tpr\tau_{(2)}\right)\nonumber\\
    &=\Omega(\omega_{(1)})\tpr\Omega(\omega_{(2)})\spr
    \Omega(\tau_{(1)})\tpr\Omega(\tau_{(2)})\nonumber\\
    &=\Omega(\omega_{(1)})\sqcup\Omega(\tau_{(1)})\tpr
    \Omega(\omega_{(2)})\sqcup\Omega(\tau_{(2)})\nonumber\\
    &=\Omega(\omega_{(1)}\tau_{(1)})\tpr\Omega(\omega_{(2)}\tau_{(2)})\nonumber\\
    &=(\Omega\tpr\Omega)\omega_{(1)}\tau_{(1)}\tpr\omega_{(2)}\tau_{(2)}\nonumber\\
    &=(\Omega\tpr\Omega)\Delta_C(\omega\tau).
  \end{align}
  Now we prove~\eqref{eq:omegarelation2} by induction.  For a forest we
  find using \eqref{eq:coprod-r} and \eqref{eq:relate}
  \begin{align*}
    \Delta_N(\Omega(\omega\tau))
    & = \Delta_N(\Omega(\omega)\sqcup\Omega(\tau))\\
    & = \Delta_N(\Omega(\omega))\spr\Delta_N(\Omega(\omega)))\\
    & = (\Omega\tpr \Omega)\Delta_C(\omega)\spr(\Omega\tpr\Omega)\Delta_C(\tau))\\
    & = (\Omega\tpr \Omega)(\Delta_C(\omega)\Delta_C(\tau))\\
    & = (\Omega\tpr \Omega)\Delta_C(\omega\tau).
  \end{align*}
  For a tree $\tau=B^+_i(\omega)\in\UT$ and $\widetilde{\omega}\in\UF$ we
  find using the definition of $\Omega$ and the recursion formulas
  \eqref{eq:delta-r} and \eqref{eq:DcRecur} that
  \begin{align*}
    \Delta_N(\Omega(\tau)) & = \Delta_N(\Omega(B^+_i(\omega))) =
    \Delta_N(B^+_i(\Omega(\omega)))\\ 
    & =  B^+_i(\Omega(\omega))\tpr\one + (\one\tpr B^+_i)\Delta_N(\Omega(\omega)) \\
    & = \Omega(\tau)\tpr \one + (\one\tpr B^+_i)(\Omega\tpr\Omega)\Delta_C(\omega)\\
    & = (\Omega\tpr\Omega) (\tau\tpr\one + (\one\tpr
    B^+_i)\Delta_C(\omega))\\
    & = (\Omega\tpr\Omega) \Delta_C(\tau).
  \end{align*}
  The proof extends to a general $a\in N$ by linearity.  {\hfill$\Box$
    \vskip10pt}
\end{proof}

A consequence of the above theorem is a connection between the
antipodes of the corresponding Hopf algebras.
\begin{corollary}
  Given $\Omega:C\rightarrow N$ is a bialgebra homomorphism then for
  $\omega\in\UF$, it follows that
  \begin{equation*}
    S_N(\Omega(\omega))=\Omega(S_C(\omega)).
  \end{equation*}
\end{corollary}
Note that symmetrization operator $\Omega$ is invertible, so
expressions for the product, coproduct and antipode of the Hopf
algebra of unordered forests can be directly expressed in terms of the
corresponding functions in the Hopf algebra of ordered forests, they
are
\begin{align*}
  \mu_C(\omega_1\tpr\omega_2)
  &=\Omega^{-1}(\mu_N(\Omega(\omega_1)\tpr\Omega(\omega_2))),\\
  \Delta_C(\omega)&=(\Omega\tpr\Omega)^{-1}\Delta_N(\Omega(\omega)),\\
  S_C(\omega)&=\Omega^{-1}(S_N(\Omega(\omega))).
\end{align*}

In the final part of this section we will elaborate on the connections
between the LS-series built on ordered trees, and their commutative
counterpart the S-series.  These series belong respectively to the
dual spaces $\Nd$ and $\Cd$, are naturally associated through the dual
map $\signature{\Omega^*}{\Nd}{\Cd}$ taking the series of ordered
forests to unordered forests.  If $\alpha\in \Nd$ and
$\beta=\Omega^*(\alpha)\in \Cd$, we find
\begin{equation*}
  \beta(\omega) = \langle \Omega^*(\alpha),\omega\rangle = 
  \langle \alpha,\Omega(\omega)\rangle = 
  \sigma(\omega) \sum_{\widetilde{\omega}\sim \omega} \alpha(\widetilde{\omega}).
\end{equation*}
On a manifold with a commutative Lie group action the elementary
differential operators $\F(\omega)$ do not depend on the
ordering. Thus we find that the S-series of $\beta$ as defined
in~\cite{murua99fsaI} equals the LS-series of $\alpha$ as given
in~\eqref{eq:lbseries},
\begin{equation*}
  \LS(\alpha) =
  \sum_{\omega\in\OF}h^{|\omega|}\alpha(\omega)\F(\omega) = 
  \sum_{\omega\in\UF}\frac{h^{|\omega|}\beta(\omega)}{\sigma(\omega)}\F(\omega) 
  = S(\beta).
\end{equation*}
This shows that the normalization $1/\sigma(\omega)$ in the
commutative case is compatible with our normalization in the
non-commutative case.

It is interesting to characterize the image of the \emph{logarithmic}
and \emph{exponential} series under $\Omega^*$.  If $\alpha$ is
logarithmic (Lemma~\ref{lem:log_exp}) then
\begin{align*}
  \beta(\one) &= \alpha(\one) = 0,\\
  \beta(\omega\tau) &= \langle\alpha,\Omega(\omega\tau)\rangle =
  \langle\alpha,\Omega(\omega)\sqcup\Omega(\tau)\rangle = 0,
  \quad\mbox{for $\omega,\tau\neq \one$},
\end{align*}
thus $\beta$ is non-zero only on trees.  If $\alpha$ is exponential
then
\begin{align*}
  \beta(\one) &= \alpha(\one) = 1,\\
  \beta(\tau_1\cdots\tau_k) 
  &= \langle\alpha,\Omega(\tau_1\cdots\tau_k)\rangle = 
  \langle\alpha,\Omega(\tau_1)\sqcup\cdots\sqcup\Omega(\tau_k)\rangle \\
  &= \alpha(\Omega(\tau_1))\cdots \alpha(\Omega(\tau_k)) =
  \beta(\tau_1)\cdots\beta(\tau_k).
\end{align*}
This is a well known condition in the composition of B-series, see
\cite{chartier05aaa}.

% ===========================================================================
\subsection{The Free Associative Algebra and frozen vector fields}
% ===========================================================================

In the proof of Proposition~\ref{prop:antipode_nr}, we defined the
Hopf algebraic structure $\Hf$ as the Free Associative Algebra built
of words over an alphabet $A$, where $\mu_F$ is the shuffle and
$\Delta_F$ the dual of the concatenation product. We will briefly
comment upon the connection between $\Hn$ and $\Hf$ in the context
of Lie group integrators.

In the theory of Lie group integrators, it is common to call constant
sections $(\M\rightarrow\g)$ \emph{frozen vector fields}. If $f$ is
frozen, then the Lie derivative $g[f] = 0$ for all vector fields $g$.
On the algebraic side, a tree $\tau\in \OT$ represents a frozen
vector field if the left grafting of anything non-constant to the tree
is zero,
\begin{equation*}
  \widetilde{\tau}[\tau] = 0, 
  \qquad\mbox{for all $\widetilde{\tau}\in \OT\backslash\{\one\}$}. 
\end{equation*}
In this case we see that the Grossman--Larson product becomes just the
concatenation product $\widetilde{\tau}\opr \tau = \widetilde{\tau}
\tau$. The freezing of certain vector fields can be understood as the
quotient $\Hf = \Hn/G$, where $G$ is the linear span of any
$\widetilde{\tau}\in \OT\backslash\{\one\}$ grafted to a frozen vector
field. As a special example, we consider the case where all
single-node trees are frozen, so that taller trees cannot be
produced. Letting $A = \{B_i^+(\one)\}_{i\in \mathcal{I}} =
\{a_i\}_{i\in \mathcal{I}}$ be the alphabet of all single-node trees,
we find from~\eqref{eq:delta-r} the following well known recursion for
the coproduct
\begin{align*}
  \Delta_F(\one) & = \one\tpr\one, \\
  \Delta_F(\omega a_i) & = \omega a_i\tpr\one+\Delta_F(\omega)\tpr
  (\one\tpr a_i).
\end{align*}
  
% ===========================================================================
\section{Concluding remarks}
% ===========================================================================

We have in this paper investigated the algebraic structure of the Hopf
algebra underlying numerical Lie group integrators. We have developed
both recursive and direct formulae for the coproduct and the antipode,
and we have in particular emphasized the connection to the Hopf
algebra of classical Butcher theory and to the free associative Hopf
algebra.  We believe that this work is of particular interest for the
construction of symbolic software packages dealing with computations
involving algebras of non-commutative derivations. The algebraic
structure of $\Hn$ is of a universal nature and should be of interest
also outside the field of numerical integration, for example, in the
renormalization of quantum field theory and the Chen--Fliess theory
for optimal control.

\begin{table}[ht]
  \centering
  \begin{equation*}
    \begin{array}{c@{\,\,}|@{\quad}l}
      \hline \\
      \omega & \Delta_N(\omega)  \\[2mm] 
      \hline \\[-1mm]
      \one & \one\tpr\one \\[1mm] 
      \ab & \ab\tpr\one+\one\tpr\ab \\[2mm]     
      \aabb & \aabb\tpr\one+\ab\tpr\ab+\one\tpr\aabb  \\[2mm]
      \ab\ab & \ab\ab\tpr\one+\ab\tpr\ab+\one\tpr\ab\ab \\[2mm]
      \aaabbb &\aaabbb\tpr\one+\ab\tpr\aabb+\aabb\tpr\ab+
      \one\tpr\aaabbb \\[2.5mm]
      \aababb & \aababb\tpr\one+\ab\ab\tpr\ab+\ab\tpr\aabb+
      \one\tpr\aababb \\[2.5mm]
      \ab\aabb & \ab\aabb\tpr\one+2\ab\ab\tpr\ab
      +\ab\tpr\aabb+\ab\tpr\ab\ab+
      \one\tpr\ab\aabb \\[2.5mm]
      \aabb\ab & \aabb\ab\tpr\one+\aabb\tpr\ab+
      \ab\tpr\ab\ab+\one\tpr\aabb\ab \\[2mm]
      \ab\ab\ab & \ab\ab\ab\tpr\one+\ab\ab\tpr
      \ab+\ab\tpr\ab\ab+\one\tpr\ab\ab\ab \\[2mm]
      \aaaabbbb & \aaaabbbb\tpr\one+\aaabbb\tpr\ab+\aabb\tpr\aabb+
      \ab\tpr\aaabbb+\one\tpr\aaaabbbb \\[2.5mm]
      \aaababbb & \aaababbb\tpr\one+\aababb\tpr\ab+\ab\ab\tpr\aabb+
      \ab\tpr\aaabbb+\one\tpr\aaababbb \\[2.5mm]
      \aabaabbb & \aabaabbb\tpr\one+\ab\aabb\tpr\ab+2\ab\ab\tpr
      \aabb+\ab\tpr\aaabbb+\ab\tpr\aababb+\one\tpr\aabaabbb \\[2.5mm]
      \aaabbabb & \aaabbabb\tpr\one+\aabb\ab\tpr\ab+\aabb\tpr
      \aabb+\ab\tpr\aababb+\one\tpr\aaabbabb \\[2.5mm]
      \aabababb & \aabababb\tpr\one+\ab\ab\ab\tpr\ab+\ab\ab\tpr
      \aabb+\ab\tpr\aababb+\one\tpr\aabababb \\[2.5mm]
      \ab\aaabbb & \ab\aaabbb\tpr\one+\ab\aabb\tpr\ab+\aabb\ab\tpr\ab+\aabb\tpr
      \ab\ab+2\ab\ab\tpr\aabb+\ab\tpr\ab\aabb+
      \ab\tpr\aaabbb+\one\tpr\ab\aaabbb \\[2.5mm]
      \aaabbb\ab & \aaabbb\ab\tpr\one+\aaabbb\tpr\ab+\aabb\tpr
      \ab\ab+\ab\tpr\aabb\ab+\one\tpr\aaabbb\ab \\[2.5mm]
      \ab\aababb & \ab\aababb\tpr\one+3\ab\ab\ab\tpr\ab+
      \ab\ab\tpr\ab\ab+2\ab\ab\tpr\aabb+
      \ab\tpr\ab\aabb+\ab\tpr\aababb+
      \one\tpr\ab\aababb \\[2.5mm]
      \aababb\ab &  \aababb\ab\tpr\one+\aababb\tpr\ab+\ab\ab\tpr
      \ab\ab+\ab\tpr\aabb\ab+\one\tpr\aababb\ab \\[2.5mm]
      \aabb\aabb & \aabb\aabb\tpr\one+\aabb\ab\tpr\ab+\ab\aabb\tpr\ab
      +\aabb\tpr\aabb+2\ab\ab\tpr\ab\ab+
      \ab\tpr\ab\aabb+\ab\tpr\aabb\ab+
      \one\tpr\aabb\aabb \\[2.5mm]
      \ab\ab\aabb & \ab\ab\aabb\tpr\one+3\ab\ab\ab\tpr\ab+
      2\ab\ab\tpr\ab\ab+\ab\ab\tpr\aabb+
      \ab\tpr\ab\ab\ab+\ab\tpr\ab\aabb+
      \one\tpr\ab\ab\aabb \\[2.5mm]
      \ab\aabb\ab & \ab\aabb\ab\tpr\one+\ab\aabb\tpr\ab+
      2\ab\ab\tpr\ab\ab+\ab\tpr\ab\ab\ab+
      \ab\tpr\aabb\ab+\one\tpr\ab\aabb\ab \\[2.5mm]
      \aabb\ab\ab & \aabb\ab\ab\tpr\one+\aabb\ab\tpr\ab+
      \aabb\tpr\ab\ab+\ab\tpr\ab\ab\ab+
      \one\tpr\aabb\ab\ab \\[2mm]
      \ab\ab\ab\ab & \ab\ab\ab\ab\tpr\one+\ab\ab\ab\tpr\ab+
      \ab\ab\tpr\ab\ab+\ab\tpr\ab\ab\ab+
      \one\tpr\ab\ab\ab\ab \\[2mm] \hline
    \end{array}
  \end{equation*}
  \caption{The coproduct in $\Hn$ for all forest up to and
    including order four. This is the dual of the product in the
    Grossman--Larson Hopf algebra. \label{tab:coproduct}}
\end{table}

\begin{table}[ht]
  \centering
  \begin{equation*}
    \begin{array}{c@{\,\,}|@{\quad}l}
      \hline \\
      \omega & S_N(\omega)  \\[2mm] 
      \hline \\[-1mm]
      \one & \one \\ 
      \ab & -\ab \\[2mm]     
      \aabb & -\aabb+2\ab\ab  \\[1.5mm]
      \ab\ab & \ab\ab \\[2mm]
      \aaabbb & -\aaabbb+2\ab\aabb+2\aabb\ab-6\ab\ab\ab  \\[2.5mm]
      \aababb & -\aababb+\ab\aabb+\aabb\ab-3\ab\ab\ab \\[2.5mm]
      \ab\aabb & \aabb\ab-3\ab\ab\ab \\[2.5mm]
      \aabb\ab & \ab\aabb-3\ab\ab\ab \\[1.5mm]
      \ab\ab\ab & -\ab\ab\ab \\[2mm]
      \aaaabbbb & -\aaaabbbb+2\ab\aaabbb+2\aaabbb\ab+2\aabb\aabb-6\ab\ab\aabb-6\ab\aabb\ab-6\aabb\ab\ab+24\ab\ab\ab\ab \\[2.5mm]
      \aaababbb & -\aaababbb+\ab\aaabbb+\aaabbb\ab+\ab\aababb+\aababb\ab-3\ab\ab\aabb-3\ab\aabb\ab-3\aabb\ab\ab+12\ab\ab\ab\ab \\[2.5mm]
      \aabaabbb & -\aabaabbb+\ab\aaabbb+\aaabbb\ab+\ab\aababb+\aababb\ab-2\ab\ab\aabb-3\ab\aabb\ab-4\aabb\ab\ab+12\ab\ab\ab\ab \\[2.5mm]
      \aaabbabb & -\aaabbabb+\ab\aababb+\aababb\ab+2\aabb\aabb-4\ab\ab\aabb-3\ab\aabb\ab-2\aabb\ab\ab+12\ab\ab\ab\ab \\[2.5mm]
      \aabababb & -\aabababb+\ab\aababb+\aababb\ab-\ab\ab\aabb-\ab\aabb\ab-\aabb\ab\ab+4\ab\ab\ab\ab \\[2.5mm]
      \ab\aaabbb & \aaabbb\ab-\ab\ab\aabb-2\ab\aabb\ab-3\aabb\ab\ab+12\ab\ab\ab\ab \\[2.5mm]
      \aaabbb\ab & \ab\aaabbb-3\ab\ab\aabb-2\ab\aabb\ab-\aabb\ab\ab+12\ab\ab\ab\ab \\[2.5mm]
      \ab\aababb & \aababb\ab-\ab\aabb\ab-2\aabb\ab\ab+6\ab\ab\ab\ab\\[2.5mm]
      \aababb\ab & \ab\aababb-2\ab\ab\aabb-\ab\aabb\ab+6\ab\ab\ab\ab \\[2.5mm]
      \aabb\aabb & \aabb\aabb-2\ab\ab\aabb-2\ab\aabb\ab-2\aabb\ab\ab+12\ab\ab\ab\ab \\[2.5mm]
      \ab\ab\aabb & -\aabb\ab\ab+4\ab\ab\ab\ab \\[2.5mm]
      \ab\aabb\ab & -\ab\aabb\ab+4\ab\ab\ab\ab \\[2.5mm]
      \aabb\ab\ab & -\ab\ab\aabb+4\ab\ab\ab\ab \\[2mm]
      \ab\ab\ab\ab & \ab\ab\ab\ab \\[2mm] \hline 
    \end{array}
  \end{equation*}
  \caption{The antipode in $\Hn$ for all forests up to and including
    order four. This is the dual of the Grossman--Larson antipode.  
    \label{tab:antipode}} 
\end{table}

\end{document}